\RequirePackage{fix-cm}
\documentclass[smallextended]{svjour3}       
\smartqed  
\RequirePackage{graphicx}
\RequirePackage{amsfonts}
\RequirePackage{amsmath}
\RequirePackage{amssymb}
\RequirePackage{mathtools}
\RequirePackage{mathabx}
\RequirePackage{dsfont}

\newcommand\tab[1][1cm]{\hspace*{#1}}

\DeclarePairedDelimiter\floor{\lfloor}{\rfloor}

\begin{document}

\title{p-composites of the main sequence of odd numbers as building blocks for $\pi(x)$}



\author{
Yuri Heymann\
}


\institute{Yuri Heymann \at
              \email{y.heymann@yahoo.com}             \\
            \emph{Present address:} 3 rue Chandieu, 1202 Geneva, Switzerland 
 }

\date{Received: date / Accepted: date}

\maketitle

\begin{abstract}
The prime-counting function $\pi(x)$ which returns the number of primes smaller or equal to a given number is a topic of interest in number theory. An algorithm based on a cyclic group isomorphic to $Z/nZ$, the so-called $Z$-functions, was proposed in view to outperform its pieers. The approach suggests a time complexity $\mathcal{O}(x^{1/2})$ in agreement with optimality of a 2-D squared adaptive-recursive algorithm. The present work is a presentation of various approaches as ascending factorization, the main sequence of odd numbers and partial sequences, T-series, counting function of prime composites, $Z$-modular forms and combinatorial aspects.

\end{abstract}

\section{Introduction}

\noindent
The prime-number theorem describing the asymptotic approximation of the distribution of primes, as defined by the prime-counting function $\pi(x)$ for the number of primes smaller or equal to $x$, was proven by a variety of authors as \cite{Poussin,Hadamar1886b}. While methods based on the prime-number theorem and Riemann hypothesis \cite{Mangoldt,Riemann1859,Hans} are of interest for computation of the prime-counting function, more common approaches are enumerated further down as the Meissel-Lehmer algorithm \cite{Meissel} and its variants and enhancements, e.g. Lehmer \cite{Lehmer}, Lagarias and Odlyzko \cite{Lagarias}, Del{é}glise and Rivat \cite{Deleglise}, etc.

The Riemann’s explicit formula is a key reference; see the Riemann-von Mangoldt's formula \cite{Mangoldt} also referred in \cite{Dittrich,Franke,Platt}, validated a while ago. ``The method was successfully validated for $\pi(x)$ calculations up to $1000$ with the first 29 pairs of complex zeros only'' as reported by Hans et al. \cite{Hans}.
 
\noindent
Other non-standard methods for calculating $\pi(x)$ are common. The Sieve of Eratosthenes is a classical example obtained by crossing out odd numbers from a list with a basis set of generators \cite{Fisher}. The Sieve of Sundaram \cite{Aiyar} and Sieve of Atkin \cite{Atkin} known for its computational hedge are variations of the former. Sharing a similitude with the aforementioned Sieves, the approach for $\pi(x)$ calculation in the present manuscript relies on the sequence of odd numbers and subsequences. Provided the sequence contains all primes larger than $a$, the number of primes belonging to $]a, x]$ is the simple difference betwen the number of odds and non-primes where the last element of the sequence is the largest smaller or equal to $x$.

A variety of methods under consideration thereoff include ascending factorization, the main sequences of odd numbers and partial sequences, T-series, counting of prime composites, modular $Z$-functions as bricks to build a better $\pi(x)$ and related topics in number theory and combinatorials.






\section{The sequence $\{S_n\}$ a key concept for primes; definitions and propositions}

The below is a proof of concept for primality, prime-counting and prime-number generation. Let us introduce the sequence $\{S_n\}$ where $n \in \mathbb{N}$, defined as follows:

\begin{equation}
S_{n+2} = S_n +6 \,,
\end{equation}

\noindent
with initial conditions

\begin{equation}
S_0 = 5 \,,
\end{equation}

\noindent
and

\begin{equation}
S_1 = 7 \,,
\end{equation}

\noindent
\textbf{Proposition 1} Any prime number larger than seven can be expressed as a number of the sequence $\{S_n\}$ added to six.

\noindent
\textbf{Proposition 2} Any prime number added to six is a prime number or the product of numbers of the sequence $\{S_n\}$ where a number can drawn a multiple amount of time. (law of composition)

\vspace{2mm}
\noindent
\underline{The law of composition, odd numbers and axiom of completness}

\vspace{1mm}
\noindent
The above two propositions stem from  the \textbf{law of composition}, i.e. the set of elements of the sequence $\{S_n\}$ including the identity element is spanned by the semigroup $G$, defined by the set $E$ which has for basis the prime numbers larger than three and neutral element 1 provided with an internal composition law which for any two elements $(x,y)$ of $E \times E$ has for result the product $x \star y \in E$.

\noindent
As a rationale any odd number which is not prime can be expressed as the product of two or more odd numbers. Given $\{U_n\}$ the sequence of odd numbers containing all primes larger than two, any element of $\{U_n\}$ is either a prime number or a composite made of the product of elements of the sequence, and where a number can be drawn a multiple amount of times. 

\vspace{1mm}
\noindent
We can show that the above sequence $\{S_n\}$ contains all primes larger than three and some composites as we can write $S_{n+1}=S_n+\delta_n$, where $S_0=5$ and $\delta_n=2$ if $n$ is even and 4 if $n$ is odd. The sequence of odd numbers $U_{n+1}=U_n+2$ where $U_0=5$ contains all primes larger than three, two being the only even prime larger than one (axiom of completness). The set of elements of the sequence $\{S_n\}$ is spanned by $\{U_n\}$ excluding points divisible by 3, e.g. $V_n= 2 \, (4 + 3 \, n) +1$ where $n= 0,1,2..$. We get $V_n=9+6 \, n$ where $n \in \mathbb{N}$, which is a sub-sequence of the odd numbers divisible by 3.

\vspace{2mm}
\noindent
\underline{The main sequence of odd numbers, partial sequences and T-series}

\vspace{1mm}

\noindent
The combinatorial approach to prime counting is based on the main sequence of odd numbers $\{U_n\}$ where $n \in \mathbb{N}$, defined as  $U_{n +1}= U_{n} + 2$ with $S_0 = 3$. As an arithmetic sequence of common difference $q=2$, the main sequences can be rewritten as $U_n= 3 + 2 \, n$ where $n \in \mathbb{N}$. 

\vspace{1mm}

\noindent
The primality sequence as a part of $\mathbb{Z}$: the sequence of odd numbers $U_{n+1}=U_{n}+2$ where $S_0 = 3$ (belonging to $\mathbb{Z}$) contains all primes larger than two (proof as above). We can also say this statement is the complementary of two being the only even prime larger than one (by set complementarity), a brick of axiomatic.

 \vspace{1mm}
\noindent
Partial sequences of odd numbers are sub-sequences of $\{U_n\}$ resulting from crossing out of odds from the main sequence, e.g. odds divisible by three, five and so on. The interest in such decomposition is to come up with simple rules for the counting of prime composites, as a brick for the prime-counting function. 

\vspace{2mm}
\noindent
The sequence of odd numbers non-divisible by three is a T-series expressed as $T_{n +2}= T_{n} + 6$, where $n \in \mathbb{N}$ and where the initial values are given by $T_0 = 5$ and $T_1 = 7$. The increment by 6 is made of the product of 2 by 3.

\vspace{2mm}

\noindent
The same idea can be extended to generalised $\mathcal{T_D}$ sequences of $\mathcal{D}$ non-divisibility, where $\mathcal{D}$ is a subset of the main sequence. As a proof of concept, see below:

\vspace{2mm}

\noindent\fbox{%
    \parbox{\textwidth}{%
The sequence of odd numbers non-divisible by five expressed as $T_{n + 3}= T_{n} + 10$, where $n \in \mathbb{N}$ and initial values $T_0 = 7$, $T_1 = 11$, $T_2=13$ and $T_4 = 19$. The increment by 10 in that $T$-series as the product of 2 by 5. 

\vspace{3mm}
\noindent
 Same idea, the sequence of odd numbers non-divisible by 3 and 5 as given by  $T_{n + 8}= T_{n} + 30$, where $n \in \mathbb{N}$ and the initial values: $T_0 = 7$, $T_1=11$, $T_2=13$, $T_3 = 17$, $T_4 = 19$, $T_5 = 23$, $T_6 = 29$ and $T_7 = 31$. The increment by 30 made of the product of 2 by 3 and 5. (valid upon..)
}
}


\vspace{3mm}
\noindent
The counting of the number of elements of the main sequence of odd numbers and relative index function as bricks for the prime-counting function are given below. 

\vspace{1mm}
\noindent
Say the function $\sigma\colon {\mathbb{ R}} \to {\mathbb{ N}}$ returns the largest element of the sequence  $\{U_n\}$ where $n=0,1,2,...$ which is smaller or equal to a real number $x$. For the function $\sigma$ to be defined, its argument must be larger or equal to the first element of the sequence. 

\vspace{2mm}
\noindent
Say $M_n$  is the number of elements of the sequence  $\{U_n\}$, $n=0,1,2...$ where the last element is smaller or equal to $x$. For the sequence of odd numbers $\{U_n\}$ defined by $S_n=3+ 2 \, n$, where $n \in \mathbb{N}$, we have:

\begin{equation}
M_n = \frac{\sigma(x)-3}{2}+1 \,.
\end{equation}

\noindent
The index function $\eta\colon {\mathbb{N}} \to {\mathbb{N}}$ which takes as argument an element of the sequence $\{U_n\}$ where $n=0,1,2...$ and returns that index is as follows:

\begin{equation}
\eta(U_n) = \frac{U_n-3}{2}  \,.
\end{equation}

\section{Ascending factorization and counting of primes composites}

\noindent
Say the $k*l$ composites are defined as the products of a given element $k$ of the main sequence $\{U_n\}$ with other elements $l$ of the sequence such that $l \geq k$. For example, the sequence of $k*l$ composites associated to $U_0 = 3$ is spanning the set: $3 \times 3$, $3 \times 5$, $3 \times 7$,  $3 \times 9$, $3 \times 11$, $3 \times 13$, etc. The sequence of $k*l$ composites associated with $U_1=5$ is spanning: $5 \times 5$, $5 \times 7$, $5 \times 9$, $5 \times 11$, etc. 
\vspace{3mm}

\noindent
By square rule composition, $k*k*l$ composites are defined as the products of a square of the main sequence $k^2$ with other elements $l$ of the sequence such that $l \geq k$.
\vspace{3mm}

\noindent
The idea underpinning ascending factorization is the decomposition of natural numbers into a product of primes ordered into ascending order, and where each prime has a corresponding multiplicity. The benefit sought by having primes sorted in ascending order is to avoid double counting resulting from a prime permutation. Say the $k_1*k_2*k_3$ composites are defined as the products of three distinct primes of the main sequence denoted $k_1$, $k_2$ and $k_3$ such that $k_1 < k_2 < k_3$. We note that any odd number can be decomposed into a product of primes of $k$-multiplicity in ascending order. 

\vspace{3mm}
\noindent
The counting function of composites denoted $\text{count}(\kappa,n)$, returns the number of times a given composite $\kappa$ occurs in the main sequence $\{U_n\}$ up to the lastest element $U_n$, $n \in \mathbb{N}$ of the sequence. 

\vspace{3mm}
\noindent
The number of non-primes of the main sequence $\{U_n\}$ or any subsequence containing all primes larger than $a$ and where the last element is defined by $U_n=\sigma(x)$  where $n \in \mathbb{N}$, can take various forms as in the below $W_n$ expression:

\begin{equation}
\begin{aligned}
W_{n} &=   \text{count}(k*l, n)  \\
& - \text{count}(k*k*l,n) + \text{count}(k*k*k,n) \\
&  - \text{count}(k*k*k*l,n) + \text{count}(k*k*k*k,n) + ....  \\ 
&  - \text{count}(k_1*k_2*l,n)  \\
& -2 \, \text{count}(k_1*k_2*k_3,n) - 3 \,  \text{count}(k_1*k_2*k_3*k_4,n)  - ... \,.
\end{aligned}
\end{equation}

\noindent
where $\text{count}(\kappa,n)$ is the number of time a given composite $\kappa$ appears in the sequence $\{U_n\}$ up to the last element $U_n$, $n \in \mathbb{N}$.

\vspace{3mm}
\noindent
Depending upon the size of the collection of odd numbers and primes multiplicity, the combinatorial arrangements of the prime composites lead to a variety of $W_n$ expressions. These variations belonging to $\Omega_n$ a topological space characterised by subspaces of the same family as a result of branching and merging of modular forms, primes multiplicity, and index $n$ of a combinatorial nature. Notwithstanding, prime factors of multiplicity higher than one, require special attention due to crossover with other factors involving potential double counting, a property of the $\Omega_n$ family, characterised by some complexity see P/P-NP problems, etc.


\vspace{3mm}
\noindent
Of a modular, arithmetic and congruent nature, $\Omega_n$ complexity can be linked with dimensionalities of space and degrees of freedom, see Venn diagram in 2-D space with n-sets, partitioning, joint functions and branching as in binary search trees, etc.


\vspace{2mm}

\noindent
The above counting rules as applied to the main sequence of odd numbers, lead to the prime-counting function $\pi(x)$ expressed as follows:

\begin{equation}
\pi(x)=M_n-W_n +m\,,
\end{equation}

\noindent
where $x> S_0$ and $n=\eta(\sigma(x))$; $m$ a natural number and $W_n$ some expression as seen above.

\vspace{3mm}

\noindent
Some examples of p-counting functions of the main sequence of odd numbers in use in ascending factorization are made available below.

\vspace{3mm}
\noindent
The counting function of the $k*l$ composites in the sequence $\{U_n\}$  is expressed as follows:

\begin{equation}
\begin{aligned}
\text{count}(k*l,n) & = \sum_{k=3, \text{k is odd}}^{\chi_n}  \floor*{\frac{S_n-k^2}{2 \, k}}+1  \\
& = \frac{\chi_n-3}{2} +1 + \sum_{k=3, \text{k is odd}}^{\chi_n}  \floor*{\frac{S_n-k^2}{2 \, k}} \,.
\end{aligned}
\end{equation}

\noindent
where $\chi_n=\eta(\sigma(\sqrt{U_n}))$.

\vspace{3mm}
\noindent
The counting function of the $k*k*l$ composites in the sequence $\{U_n\}$ is expressed as follows:

\begin{equation}
\begin{aligned}
\text{count}(k*k*l,n) & = \sum_{k=3, \text{k is odd}}^{\chi_n}  \floor*{\frac{n}{k^2}+ \frac{\xi(k)}{k^2}} \,,
\end{aligned}
\end{equation}

\vspace{3mm}
\noindent
where $\chi_n=\eta(\sigma(\sqrt{U_n}))$ and the function $\xi \colon \mathbb{N} \rightarrow \mathbb{R}$ is defined as follows:

\begin{equation}
\xi(k) =\frac{3-k^2}{2}\,.
\end{equation}

\vspace{3mm}
\noindent
The counting function of the $k^j$ composites in the sequence $\{U_n\}$ where $j \in \, \mathbb{N}^{*}$ is expressed as follows:

\begin{equation}
\begin{aligned}
\text{count}(k^j,n) =  \floor*{\frac{{U_n}^{1/j}-3}{2}} +1 \,,
\end{aligned}
\end{equation}

\vspace{3mm}
\noindent
where $\chi_n=\eta(\sigma(\sqrt{U_n}))$.

\vspace{3mm}
\noindent
The counting function of composites formed of the product of primes of multiplicity larger than one in the sequence $\{U_n\}$ (e.g. $k_1*k_2*k_3$ seen above) can be viewed as a family of elliptic curves, see B{\'e}zout's identity for reference.


\section{$Z$-functions as bricks for prime composites and factorization}

The axiom of completness tells us that all primes excluding two are contained in the main sequence of odd numbers $\{U_n\}$ where $n \in \mathbb{N}$. From internal composition, the number of primes in the sequence is the difference between the overall number of elements $M_n$ and non-primes $W_n$, defined as all possible composites built from the product of numbers of the sequence of an integer multiplicity, and where the product is smaller or equal to the last term $U_n$ of the sequence, where $U_n=\sigma(x)$. The idea is to count all composites from the sequence $\{U_n \}$, avoiding double counting as a result of non-coprime factors.

\vspace{3mm}
\noindent
Let us introduce the $3$-composites defined as the sequence of the products of the prime $3$ with the other elements of the sequence $U_k$ where $U_k \geq 3$ and $k \in \, \mathbb{N}$. The sequence is as follows: $3 \times 3$, $3 \times 5$, $3 \times 7$, $3 \times 9$, etc. We define the below function for counting of the 3-composites:

\begin{equation}
\text{count}(3,n) = \floor*{\frac{n}{3}} \,,
\end{equation}

\noindent
where $n=\eta(\sigma(x))$. 

\vspace{3mm}
\noindent
Let us say the $5$-composites are defined as products of primes $p_k \geq 5$ raised to powers $m_k \geq 1$, containing at  least a 5 of arbitrary multiplicity as the sequence of prime 5 multiplied by odds non-divisible by 3, where the sequence is as follows: $5 \times 5$, $5 \times 7$, $5 \times 11$, $5 \times 13$, $5 \times 17$, etc. This sequence starts at $5*5$ and has succesive jumps $+15$, $+15$, $+15$, etc. The modular function for the counting of the $5$-composites in $\{U_n\}$ is as follows:

\begin{equation}
\text{count}(5,n) = \mathds{1}_{n \geq 11}  \floor*{ \,\frac{1}{3} \floor*{\frac{n-11}{5}}+\frac{1}{3} } \,,
\end{equation}

\noindent
where the integer $11$ in (13) is coming from $\sigma(5^2/2)$.

\vspace{3mm}
\noindent
Let us introduce the $7$-composites defined as the sequence formed of the products of prime $7$ with the ordered sequence of odds non-divisible by 3. The sequence is as follows: $7 \times 7$, $7 \times 11$, $7 \times 13$, $7 \times 17$, etc. This sequence starts at $7*7$ and has succesive jumps $+14$, $+7$, $+14$, $+7$, etc. Its modular function for the counting of the $7$-composites in $\{U_n\}$ is as follows:

\begin{equation}
\text{count}(7,n) = \mathds{1}_{n \geq 23} \left(1+ \floor*{\frac{n-23}{7}} - \floor*{ \,\frac{1}{3} \floor*{\frac{n-23}{7}}+\frac{2}{3} }\right) \,,
\end{equation}

\noindent
where the integer $23$ in (14) is coming from $\sigma(7^2/2)$.

\vspace{3mm}
\noindent
Let us introduce the $11$-composites defined as the sequence of the products of prime $11$ with ordered sequencee of odds non-divisible by 3. The sequence is as follows: $11 \times 11$, $11 \times 13$, $11 \times 17$, etc. This sequence starts with $11*11$ and has succesive jumps $+11$, $+22$, $+11$, $+22$, etc. Thus, the counting function of the $11$-composites in $\{U_n\}$ is as follows:

\begin{equation}
\text{count}(11,n) = \mathds{1}_{n \geq 59} \left(1+ \floor*{\frac{n-59}{11}} - \floor*{ \,\frac{1}{3} \floor*{\frac{n-59}{11}}+\frac{1}{3} }\right) \,,
\end{equation}

\vspace{2mm}
\noindent
where the integer $59$ in (15) is coming from $\sigma(11^2/2)$.

\vspace{3mm}

\vspace{2mm}
\noindent
An example of recursive algorithm for the generation of the first N primes on the basis of the sequence of odd numbers and above z-modular forms, is provided as is in the below pseudo code. 

\newpage

\vspace{4mm}
\noindent
\textbf{Pseudo-code for generation of the first N primes}
\vspace{4mm}

\noindent\fbox{%
    \parbox{\textwidth}{%
    
    int getMainSequenceElement(int n):
    	return 3 + 2*n;
    	\vspace{2mm}
    	
    int getIndexAtElement(int element):
    	return (element - 3)/2;
     \vspace{2mm}
    
    int getSigma(double x):
    	return (int)((x-3.0)/2.0);
    	
    	\vspace{2mm}
    	
    arraylist$<$int$>$ getFirstNPrimes(int N) \\
     \{ \\
    	\tab arraylist$<$int$>$ primeSequence; \\
    	\tab arraylist$<$int$>$ moduloSequence;
    	
    	\vspace{2mm}
    	
    	\tab //initialize the sequences\\
    	\tab primeSequence.add(3);\\
    	\tab primeSequence.add(5);\\
    	\tab moduloSequence.add(3);\\
    	\tab moduloSequence.add(5);\\
    	
    	\tab /*  extract primes from the first partition\\
    	 \tab \hspace{1mm} * anchor points: prime 5 and prime 7\\
    	 \tab \hspace{1mm} * starts at 5 and ends at 7*7 
    	 \tab */
    	 
    	 \vspace{1mm}
    	 \tab // anchor points\\
    	 \tab int primeA = 5;\\
    	 \tab int primeB = 7;\\
    	
    	    	 \vspace{1mm}
    	 \tab int startIndex = getIndexAtElement(primeA*primeA);\\
    	 \tab int endpoint; // last index of the partition\\
    	 
    	 \tab int index = startIndex;\\
        \tab int indexForPartitions = 1;\\
        \tab int lastElement = 1; \\
        
            \tab while (primeSequence.size() $<=$ N)\\
       \tab  \{  \\
           \tab \tab  if (indexForPartitions $>$ 1)\\
            \tab \tab \{ \\
            	\tab \tab \tab //index from last element of prior partition\\
                \tab \tab \tab index = getIndexAtElement(lastElement * primeA); \\
            \tab \tab \}  \\ \\
            \tab \tab endpoint = getIndexAtElement(primeA * primeB * primeB);\\      
            
              \tab \tab  while (index $<$ endpoint - primeA)\\
            \tab \tab \{                    \\         
                
    	... continue
    	}
}

\newpage

    	\noindent\fbox{%
    \parbox{\textwidth}{%
    	    	             \vspace{5mm}
cont..... \\    
               \tab \tab \tab int nextPrime;\\
               \tab \tab \tab index += primeA; // periodicity            \\
                \tab \tab \tab nextPrime = getMainSequenceElement(index) / primeA;\\
               \tab \tab \tab  boolean isPrime = true; \\
    	    	  \tab \tab \tab  // pruning \\
                \tab \tab \tab  for (int module : moduloSequence)\\
                \tab \tab \tab \{ \\
                \tab \tab \tab \tab   if (nextPrime \% module == 0) \\
                 \tab \tab \tab \tab   \{ \\
                 \tab  \tab \tab \tab \tab      isPrime = false;\\
                 \tab \tab \tab \tab \tab      break;\\
                  \tab \tab \tab \tab \}  \\
                \tab \tab \tab \} \\ \\
                \tab \tab \tab if (isPrime)\\
               \tab \tab \tab  \{ \\
               \tab \tab  \tab \tab	primeSequence.add(nextPrime);\\
               \tab \tab \tab  \} \\
                \tab \tab \}  \\
                
           \tab \tab  // rollover \\
           \tab \tab moduloSequence.add(primeB);\\
            \tab \tab primeA = primeB;\\
            \tab \tab int position = primeSequence.indexOf(primeB);\\
           \tab \tab  primeB = primeSequence.get(position + 1);\
            
          \tab \tab  indexForPartitions++;\\
          \tab \tab  lastElement = primeSequence.get(primeSequence.size() - 1);\\
         \tab \} \\ \\
        \tab return primeSequence; \\
        \} 
   }
}

\bibliographystyle{spmpsci}      
\bibliography{Zeta2}   

%
%

\end{document}